# SIMULATION CASE STUDIES ON PERIOD OPTIMAL POWER FLOW

Zongjie Wang
Ge Guo
C. Lindsay Anderson

Cornell University
111 Wing Dr.
Ithaca, NY, 14850, US

## ABSTRACT

Traditional optimal power flow (OPF) describes the system performance only in a single snapshot while the resulting decisions are applied to an entire time period. Therefore, how well the selected snapshot can represent the entire time period is crucial for the effectiveness of OPF. In this paper, period OPF is proposed with a given time period that partitioned into consecutive linear time intervals, in which bus voltages and branch powers are linear mappings of time. The optimality and operational limits in each time interval are precisely represented by its median and two end snapshots, through which the system performance of OPF model is significantly improved. Simulation case studies on a modified IEEE 118-bus system have demonstrated the simplicity and effectiveness of the proposed model.

## 1 INTRODUCTION

The conception of optimal power flow (OPF) was proposed by J. Carpentier in 1962 and has brought huge impact to the power system field in the world (Carpentier 1962). Different from heuristic rules for power generation scheduling such as equal rate of incremental (Ringlee and Williams 1962), the number of inequality constraints in the OPF model may be enormous, which makes OPF a large-scale nonlinear optimization problem.

There have been unwavering efforts on reliable and effective solution methods for the OPF problem (Liang et al. 2011; Souza et al. 2017). Take the Newton method (Sun et al. 1984) and the interior-point method (Capitanescu and Wehenkel 2013) that deals with inequality constraints as two typical examples, OPF methods have become more developed before the end of last century. By satisfying branch power constraints and voltage constraints, OPF has guaranteed both the systems' security and the voltage quality. Therefore, OPF has become an essential advanced tool in Energy Management System (EMS) (Maria and Findlay 1987), and has been widely applied into power systems' dispatch center.

Generation scheduling made by system operators are based on the forecasts from loads and intermittent power sources (IPSs) (Wang and Guo 2018). These forecasts are usually predicted on discrete snapshots, while output levels between two consecutive discrete snapshots are approximated by linear interpolation. Thus the power output is a piecewise linear function of time from the perspective of the continuous temporal domain.

In short-term (day-ahead) scheduling, based on power forecast, system operators solve the OPF problem on a particular snapshot, yielding a short-term optimal generation schedule that periodically adjusts outputs of controllable power generation assets. In real-time dispatch, based on the real-time power forecast (Reddy et al. 2015; Wang and Guo 2018), the OPF problem is sequentially solved for certain future snapshots, thus the real-time optimal generation schedule is made for a given time period. For either of them, the OPF solution is performed on discrete snapshots. Hence a fundamental characteristic OPF is proposed to specify proper discrete snapshots to represent the continuous time period (Wang and Guo 2017).

The power flow distributions are not identical not only at these specified snapshots but also at those non-specified snapshots. All the operational constraints satisfied at specified snapshots do not ensure these constraints are satisfied at other non-specified snapshots, especially for binding constraints at specified snapshots. Taking the power balancing as a preliminary, OPF should ensure power and voltage limits to be satisfied at any snapshot in any given period (Gill et al. 2014; Gerstner et al. 2016; Wang and Guo 2017; Kourounis et al. 2018). However, the rigorous OPF model is in nature a large-scale dynamic optimization problem which may be costly to solve. Hence it is crucial to wisely specify proper discrete points that ease the computational burden and effectively represent the continuous time period.

Note that, our objective is to orient to a time period, but we specify some discrete snapshots as representatives. This indeed is a tradeoff between its performance and efficiency according to the characteristics of power systems. Hence the key is to specify the optimal discrete snapshots based on the tradeoff.

To this end, we propose a consecutive linearization technique of OPF. Given a time period, we partition it into consecutive linear time intervals. Each linear time interval has an approximately linear relationship between system states and time. Therefore, it can be accurately represented by its median and two end snapshots. Namely, we extract the objective function from median snapshot and check operational limits on the two end snapshots. Since the proposed OPF is better oriented to a time period, it is also referred to as the Period Optimal Power Flow (POPF).

The remainder of this paper is organized as follows. In Section 2, we review the nodal power flow equations and derive its generalized form. Linear time intervals are defined in Section 3 wherein some relevant properties are studied. Characteristics of median and two end snapshots are discussed in Section 4. In Section 5, we present the POPF formulation based on the previous discussions. Simulation case studies are recorded in Section 6.

## 2 POWER FLOW EQUATIONS

### 2.1 Equivalent Load

IPSs (such as wind and solar power sources) have been rapidly integrated into power systems due to economic and environmental reasons (Markogiannakis 2014; Gerstner et al. 2016). Their outputs are not totally controllable thus IPSs can be regarded as negative loads. Note that loads themselves are also uncertain and usually uncontrollable. Loads cannot be easily shed unless they are active loads. Since IPSs do not generate revenue when shed, so the full use of renewable energy assets is normally assumed. Therefore, we combine loads and IPSs in this paper as equivalent loads.

Let $P_i^e(t)$, $Q_i^e(t)$ respectively denote the nodal net active and reactive power at time $t$ and node $i$, defined by the following equations:

$$\forall t, \forall i : \begin{cases} P_i^e(t) = P_i^l(t) - P_i^{ips}(t) \\ Q_i^e(t) = Q_i^l(t) - Q_i^{ips}(t) \end{cases}, \quad (1)$$

where superscripts $l$ and $ips$ denote, respectively, terms associated with loads and IPSs.

### 2.2 Nodal Power Flow Equation

Without loss of generality, any node in the power system can either be a *power source node* in which the power injection is adjustable or an *equivalent load node* in which the power injection is uncontrollable. Note that a power source node may contain loads or IPSs, but an equivalent load node typically has no controllable power source.

Let $P_i(t), Q_i(t)$ respectively denote the calculated values of the active and reactive power injections at time $t$ and node $i$. These quantities are nonlinear functions of the node voltage[1]:

$$\begin{cases} P_i(t) = f_P(\dot{V}(t)) \\ Q_i(t) = f_Q(\dot{V}(t)) \end{cases}, \quad (2)$$

where $\dot{V}(t)$ is the vector of nodal complex voltages. Here we take the rectangular coordinates and the nodal electrical potential vector is expressed as

$$\forall t : \boldsymbol{x}(t) = [e_1(t), f_1(t), \ldots, e_n(t), f_n(t)]^T, \quad (3)$$

where $n$ is the number of nodes and $e$ and $f$ are real and imaginary parts of complex voltages.

The nodal power flow equation for the entire power grid at time $t$ is further obtained:

$$\forall t : \boldsymbol{f}(\boldsymbol{x}(t)) = \boldsymbol{g}(t), \quad (4)$$

where $\boldsymbol{g}(t)$ is given by

$$\forall t : \boldsymbol{g}(t) = -\boldsymbol{S}^e(t) + \boldsymbol{P}^c(t) + \boldsymbol{V}^c(t), \quad (5)$$

where $\boldsymbol{S}^e(t)$ is the equivalent load power vector, $\boldsymbol{P}^c(t)$ and $\boldsymbol{V}^c(t)$ are the active power vector and node voltage vector for controllable power sources, respectively.

It is evident to see that the node power flow equation is nonlinear and time varying.

## 3 LINEAR TIME INTERVALS

### 3.1 Linear Time Intervals

The derivative of the nodal power flow equation with respect to time can be expressed as

$$\forall t : \boldsymbol{J}(t) \frac{d\boldsymbol{x}(t)}{dt} = \frac{d\boldsymbol{g}(t)}{dt}, \quad (6)$$

where $\boldsymbol{J}(t)$ is the Jacobian matrix

$$\forall t : \boldsymbol{J}(t) = \frac{\partial \boldsymbol{f}(\boldsymbol{Y}, \boldsymbol{x}(t))}{\partial \boldsymbol{x}(t)}. \quad (7)$$

**Proposition 1. (Temporal Increment Equation)** *If the Jacobian matrix $\boldsymbol{J}(t)$ is a continuous mapping in the time interval $T_k = [t_{k-1}, t_k]$, then the following temporal increment equation holds*

$$t_m \in T_k : \boldsymbol{J}(t_m)(\boldsymbol{x}(t_k) - \boldsymbol{x}(t_{k-1})) = \boldsymbol{g}(t_k) - \boldsymbol{g}(t_{k-1}), \quad (8)$$

where $t_m = \dfrac{t_{k-1} + t_k}{2}$ is the median snapshot of the interval $T_k$.

**Proof.** Integrating (6) gives

$$\int_{t_{k-1}}^{t_k} \boldsymbol{J}(t) d\boldsymbol{x}(t) = \int_{t_{k-1}}^{t_k} d\boldsymbol{g}(t) = \boldsymbol{g}(t_k) - \boldsymbol{g}(t_{k-1}). \quad (9)$$

Since $\boldsymbol{J}(t)$ is continuous, it follows from the integral mean value theorem that

$$\int_{t_{k-1}}^{t_k} \boldsymbol{J}(t) d\boldsymbol{x}(t) = \boldsymbol{J}(t_m)(\boldsymbol{x}(t_k) - \boldsymbol{x}(t_{k-1})). \quad (10)$$

Hence the proposition holds and this end the proof. □

Let $\boldsymbol{w}$ be a constant vector, and let $t, t_* \in T_k$ be two arbitrary time points within the time interval $T_k$. If the time interval $T_k$ satisfies the following two conditions:

(1) The given nodal vector value $\boldsymbol{g}(t)$ is a linear function of time:

$$t, t_* \in T_k : \boldsymbol{g}(t) = \boldsymbol{g}(t_*) + (t - t_*)\boldsymbol{w}. \quad (11)$$

---

[1] In this paper, we assume system topology and parameters do not change in the considered time period. When the system topology changes, the whole time period should be partitioned into smaller intervals.

(2) The Jacobian matrix $J(t)$ is approximated[2] by a constant matrix $\mathfrak{J}$ in the time interval $T_k$:
$$t \in T_k : J(t) \approx \mathfrak{J} . \tag{12}$$
then the time interval $T_k$ is a *linear time interval*.

**Proposition 2. (Power Flow Linearization)** *In a linear time interval $T_k$, the node voltage vector $x(t)$ can be approximated as a linear function of time*
$$t, t_* \in T_k : x(t) \approx x(t_*) + (t - t_*) \mathfrak{J}^{-1} w . \tag{13}$$

According to Proposition 1 and the defined linear time interval, Proposition 2 holds. □

The linear time interval indeed exists unless there is an abrupt change on node power within the time interval. A sufficiently short time interval can not only ensure the linearity between nodal power and time, but also guarantee the linearity between nodal power and voltages, that is, it also guarantees the Jacobian matrix to be approximated as constant matrix.

The following discussion shows that the assumption of a linear time interval is also practically consistent with the actual characteristics of power systems.

### 3.2 Modeling for Equivalent Load Power Forecasts

The equivalent load power forecasts are the basis for the power generation scheduling, thus multiple equally distributed snapshots of equivalent loads are generated. We connect these snapshots with the values of power injections via line segments, a piecewise linear power curve is thus generated, which approximately describes the equivalent load power curve.

The interval between adjacent forecasting timepoints (each line segment) is referred to as a forecasting time interval. Thus, an $N+1$ combined equally distributed timepoints $t_k (k = 0,1,2,...,N)$ partitions into $N$ equally distributed forecasting time interval $T_k (k = 0,1,2,...,N)$ that is expressed as
$$T_k = [t_{k-1}, t_k], \tag{14}$$
where
$$\forall T_k : t_k - t_{k-1} \equiv \Delta T = constant . \tag{15}$$

There are only two equivalent load forecasting values in each time interval $T_k$, which specifies at the two end snapshots. A linear equation for the node equivalent load power is given by
$$\forall i, t \in T_k : \begin{bmatrix} P_i^e(t) \\ Q_i^e(t) \end{bmatrix} = \begin{bmatrix} P_i^e(t_{k-1}) \\ Q_i^e(t_{k-1}) \end{bmatrix} + (t - t_{k-1}) \begin{bmatrix} w_{ik}^p \\ w_{ik}^q \end{bmatrix}, \tag{16}$$
where $w_{ik}^p$ and $w_{ik}^q$ are slopes of the active and reactive power, which are constant values with respect to forecasting time intervals and is given by
$$\forall T_k, \forall i : \begin{bmatrix} w_{ik}^p \\ w_{ik}^q \end{bmatrix} = \frac{1}{\Delta T} \begin{bmatrix} P_i^e(t_k) - P_i^e(t_{k-1}) \\ Q_i^e(t_k) - Q_i^e(t_{k-1}) \end{bmatrix} . \tag{17}$$

(16) indicates that the node equivalent load vector $S^e(t)$ is a linear function of time.

### 3.3 Temporal Modeling of Active Power for Controllable Power Sources

In the dispatch system, the power generation scheduling for a forecasting time interval is determined based on a specified snapshot therein. However, the equivalent load power may deviate from their predicted values and the power generation needs to be adjusted accordingly to achieve the power balance. In this paper, we assume the deviation of equivalent load power is small in the same forecasting time interval[3].

---

[2] The precision of this approximation affects the number of width of each linear time internal, which can be empirically specified.

[3] If the equivalent load power variation in a given forecast time interval is large, the interval can be divided into two smaller time intervals

For small variations in the equivalent load power, the frequency control hierarchy of power systems will track the frequency deviation and adjust the controllable power. For both primary and secondary frequency modulation, it is assumed that the active power adjustment from controllable power source $\Delta P_i^c$ (Trovato et al. 2017) is a linear function of the grid frequency deviation $\Delta f$, that is,

$$\forall i : \Delta P_i^c(t) = -K_i \Delta f(t) , \tag{18}$$

where $K_i$ is the frequency modulation coefficient.

There are two reasons that cause grid frequency deviations $\Delta f(t)$: Inherent power imbalance from the power forecasting error and the power imbalance from time interval. In this paper, for simplicity, only the latter is considered to highlight the temporal issue, which is denoted by $\Delta f_k(t)$.

At time $t \in T_k$, the frequency deviation $\Delta f_k(t)$ is proportional to the power imbalance $\Delta P^e(t)$:

$$t \in T_k : \Delta f_k(t) = k \Delta P^e(t) , \tag{19}$$

where

$$t \in T_k : \Delta P^e(t) = \sum_{i=1}^{n} \Delta P_i^e(t) . \tag{20}$$

Considering (16),

$$\Delta P^e(t) = (t - t_{k-1}) \sum_{i=1}^{n} w_{ik}^p . \tag{21}$$

Therefore,

$$t \in T_k : \Delta f_k(t) = (t - t_{k-1}) \sum_{i=1}^{n} w_{ik}^p . \tag{22}$$

Substituting (22) into (18) gives

$$\forall i, t \in T_k : \Delta P_i^c(t) = -K_{iK}(t - t_{k-1}) , \tag{23}$$

where

$$K_{iK} = K_i \sum_{i=1}^{n} w_{ik}^p . \tag{24}$$

(23) indicates that the active power vector $\boldsymbol{P}^c(t)$ of controllable power sources is a linear function of time over each time interval.

For each forecasting time interval, the voltage magnitude variation of the controllable power sources is very small, thus the node voltage magnitude vector $\boldsymbol{V}^c(t)$ for controllable power sources is assumed to be constant.

### 3.4 Time Interval Construction

To ensure that the time intervals are linear, the time horizon can be divided using various interval lengths, and the forecast time intervals can be integrated. This is referred to as time interval construction.

Time interval construction is performed by assuming that the Jacobian matrix is approximately constant. Such an approximation indicates that in time interval $T_k$, the variation of the equivalent load power is small thus, both conditions for a linear time interval are satisfied.

Let $\mu$ be a positive constant that is sufficiently small. For each forecasting time interval $T_k$, if the norm of the Jacobian matrix at the beginning and the end of the interval (denoted as $\|\boldsymbol{J}(t_{k-1})\|$ and $\|\boldsymbol{J}(t_k)\|$, respectively) satisfy

$$\|\boldsymbol{J}(t_{k-1})\| - \|\boldsymbol{J}(t_k)\| < \mu , \tag{25}$$

then the Jacobian matrix is considered to be constant; otherwise, the time interval $T_k$ should be bisected into two smaller ones. Such a process will be iterated until (25) is satisfied.

There are various matrix norms that could be used. The 2-norm $\|\cdot\|_2$ requires the calculation of the matrix eigenvalues, which is computationally intensive. Hence the following infinity norm $\|\cdot\|_\infty$ is applied and is readily calculated.

$$\|\boldsymbol{J}\|_\infty = \max\left\{\sum_{j=1}^n |J_{1j}|, \sum_{j=1}^n |J_{2j}|, \ldots, \sum_{j=1}^n |J_{nj}|\right\}, \tag{26}$$

## 4  MEDIAN AND TWO END SNAPSHOTS

### 4.1  Median Snapshot

To make a power generation plan, values are calculated at a particular snapshot within each linear time interval $T_k$. In this study, the median snapshot of the time interval, denoted by $t_k^m$ is considered

$$\forall T_k : t_k^m = \frac{1}{2}(t_{k-1} + t_k). \tag{27}$$

The median snapshot is selected as the basis for optimization since it satisfies energy integral property.

The nodal energy vector $\boldsymbol{E}_k$ for time interval $T_k$ is given by the integral with respect to the node active power vector

$$\boldsymbol{E}_k = \int_{t_{k-1}}^{t_k} (\boldsymbol{P}^c(t) - \boldsymbol{P}^e(t))dt. \tag{28}$$

Since $\boldsymbol{P}^c(t)$ and $\boldsymbol{P}^e(t)$ are linear functions of time, (28) can be expressed as

$$\boldsymbol{E}_k = (\boldsymbol{P}^c(t_k^m) - \boldsymbol{P}^e(t_k^m))\Delta T. \tag{29}$$

In other words, the integral of power from median snapshot over a time interval accurately represents the total energy consumption.

Note that the utility maximization in OPF is directly related to energy generations and consumptions. Therefore, the median snapshot is adopted to characterize the objective function of the OPF model.

Depending on the type of the integrand function, the integral mean is not necessarily equal to the arithmetic mean. In a linear time interval $T_k$, the power vector is a linear function of time, and therefore, the median power is the arithmetic mean of the power at the beginning and the end of the interval, as shown in Figure 1. For a linear time interval $T_k$, the distance between the two endpoint values of the power and the median value of the power are equal, which ensures a minimum of power adjustments in frequency modulation.

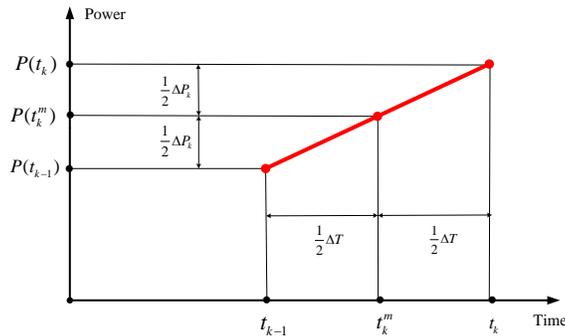

Figure 1: Equidistant property of median load pattern.

Moreover, we show in the following proposition that the voltage differences between the median and the two end snapshots are equal over a given time interval.

**Proposition 3. (Voltage Isometry Proposition)** *The voltage differences between each of the two endpoint values of the voltage vector in time interval $T_k$ and the median value of the voltage vector are approximately equal, that is,*

$$\forall T_k : \left|x(t_k) - x(t_k^m)\right| \approx \left|x(t_{k-1}) - x(t_k^m)\right|. \tag{30}$$

**Proof.** According to Proposition 2, the value at the median point is

$$x(t) \approx x(t_k^m) + (t - t_k^m)\mathfrak{J}^{-1}w, \tag{31}$$

where $w$ is a constant vector. Therefore,

$$\begin{aligned} x(t_k) - x(t_k^m) &\approx (t_k - t_k^m)\mathfrak{J}^{-1}w \\ x(t_{k-1}) - x(t_k^m) &\approx (t_{k-1} - t_k^m)\mathfrak{J}^{-1}w \end{aligned}. \tag{32}$$

Since

$$t_k - t_k^m = -(t_{k-1} - t_k^m), \tag{33}$$

it follows that

$$x(t_k) - x(t_k^m) \approx -(x(t_{k-1}) - x(t_k^m)). \tag{34}$$

### 4.2 Two End Snapshots

In a linear time interval $T_k$, nodal complex voltages are linear functions of time. Such linear functions must be either monotonic increasing or decreasing. Hence, extremal values in nodal voltage magnitudes and branch power flows occur at the two end snapshots. Therefore, the endpoint values of the function in a linear time interval $T_k$ are used to check inequality operational constraints in the OPF model. If these constraints are satisfied at the specified two end snapshots, then they are also satisfied at any other non-specified snapshots over the time interval $T_k$.

Since median snapshot has the voltage isometry property, it is sufficient to consider voltage limits on only one end snapshot.

## 5   PERIOD OPTIMAL POWER FLOW ALGORITHM

Based on the above analysis, we present in this section the scheme of the proposed POPF. Given an interested time period, we partition it into linear time intervals as described in subsection 3.4. Thereupon, the OPF is to minimize the cost of median snapshot and to check voltage and power limits for two end snapshots. Note that such a representation is rigorous as long as the linear assumption holds. The median snapshot represents the average energy consumption during linear time intervals, and the extreme values of nodal voltages and branch powers only happen at two end snapshots.

Specifically, the three-snapshot optimization based POPF is obtained as follows:

(1) Median snapshot based objective function: the objective function $f(\cdot)$ is evaluated at the median snapshot for each linear time interval $T_k$. The function could be a minimum of generation cost, the market cost, or transmission loss. It is expressed generally as

$$\forall T_k : \min f(U, X, t_k^m), \tag{35}$$

where $U$ is a vector of control variables, $X$ is a vector of state variables, and $t_k^m$ is the median snapshot.

(2) Equality constraints for median snapshot: For any linear time interval , since the inequality constraints of median snapshot necessarily satisfy if these constraints are satisfied at two endpoints, thus only the power flow equality constraints are incorporated into the median snapshot:

$$\forall T_k : h(U, X, t_k^m) = 0. \tag{36}$$

(3) Constraints for two end snapshots: for time interval $T_k$, the equality and inequality constraints must be satisfied at the two end snapshots $t_{k-1}, t_k$. The power flow equality constraints are

$$\forall T_k : \begin{cases} h(U, X, t_{k-1}) = 0 \\ h(U, X, t_k) = 0 \end{cases}. \tag{37}$$

The voltage inequality constraints are

$$\forall T_k : \begin{cases} \boldsymbol{V}_{\min} \leq \boldsymbol{V}(t_{k-1}) \leq \boldsymbol{V}_{\max} \\ \boldsymbol{V}_{\min} \leq \boldsymbol{V}(t_k) \leq \boldsymbol{V}_{\max} \end{cases}. \tag{38}$$

The controllable power inequality constraints are

$$\forall T_k : \begin{cases} \boldsymbol{P}^c_{\min} \leq \boldsymbol{P}^c(t_{k-1}) \leq \boldsymbol{P}^c_{\max} \\ \boldsymbol{Q}^c_{\min} \leq \boldsymbol{Q}^c(t_{k-1}) \leq \boldsymbol{Q}^c_{\max} \\ \boldsymbol{P}^c_{\min} \leq \boldsymbol{P}^c(t_k) \leq \boldsymbol{P}^c_{\max} \\ \boldsymbol{Q}^c_{\min} \leq \boldsymbol{Q}^c(t_k) \leq \boldsymbol{Q}^c_{\max} \end{cases}, \tag{39}$$

and the line power inequality constraints are

$$\forall T_k : \begin{cases} \boldsymbol{P}^{line}(t_{k-1}) \leq \boldsymbol{P}^{line}_{\max} \\ \boldsymbol{P}^{line}(t_k) \leq \boldsymbol{P}^{line}_{\max} \end{cases}. \tag{40}$$

The solution for the median snapshot from the previous time interval is taken as the initial value of the median snapshot for the next time interval, that is,

$$\begin{cases} \boldsymbol{P}^c_0(t^m_k) = \boldsymbol{P}^c(t^m_{k-1}) \\ \boldsymbol{V}^c_0(t^m_k) = \boldsymbol{V}^c(t^m_{k-1}) \end{cases}, (k = 1, 2, ..., N), \tag{41}$$

where $\boldsymbol{P}^c_0$ and $\boldsymbol{V}^c_0$ represent the initial values of active power vector and node voltage magnitude vector for the controllable power sources, respectively; $N$ represents the total number of end snapshots. Therefore, the proposed POPF that is oriented towards a given time period is finally obtained.

## 6    NUMERICAL SIMULATIONS

In this section, power systems with high proportional IPSs is assumed. In the power balancing process, the consecutive connection (Wang and Guo 2018) of short-term scheduling and real-time dispatch are taken into account. The optimal power generation is solved by the proposed POPF based on three representative snapshots. Our results are compared to those from a traditional optimal power flow (TOPF) algorithm which considers one single snapshot in each time interval.

For each time interval, the POPF algorithm optimizes the solution at the median snapshot subject to constraints at two end snapshots. It should be noted that the TOPF algorithm has only one computational snapshot per time interval. Their comparison is shown in Table 1. The objective function of the optimization was to minimize the cost of power generation.

Table 1: Comparison of two optimal power flow models.

| Model | Objective Function | Constraints | Time scan number(s) |
|---|---|---|---|
| POPF | Minimize generation cost on the medium snapshot | Constraints of first and two end snapshots | 3 |
| TOPF | Minimize generation cost on a single snapshot | Constraints of single snapshot | 1 |

### 6.1    System Designing

#### 6.1.1  Simulation System

The modified IEEE 118-bus system is simulated, where 20 wind generators and 10 solar generators are added. Assume that there is no curtailment of renewable generation, then the total power generation from IPSs is 30% of the total demand, with the wind power generation of being 60% of the total IPS power output. Thus a power system with high proportional IPSs is presented. The diagram of the simulation system is shown in Figure 2.

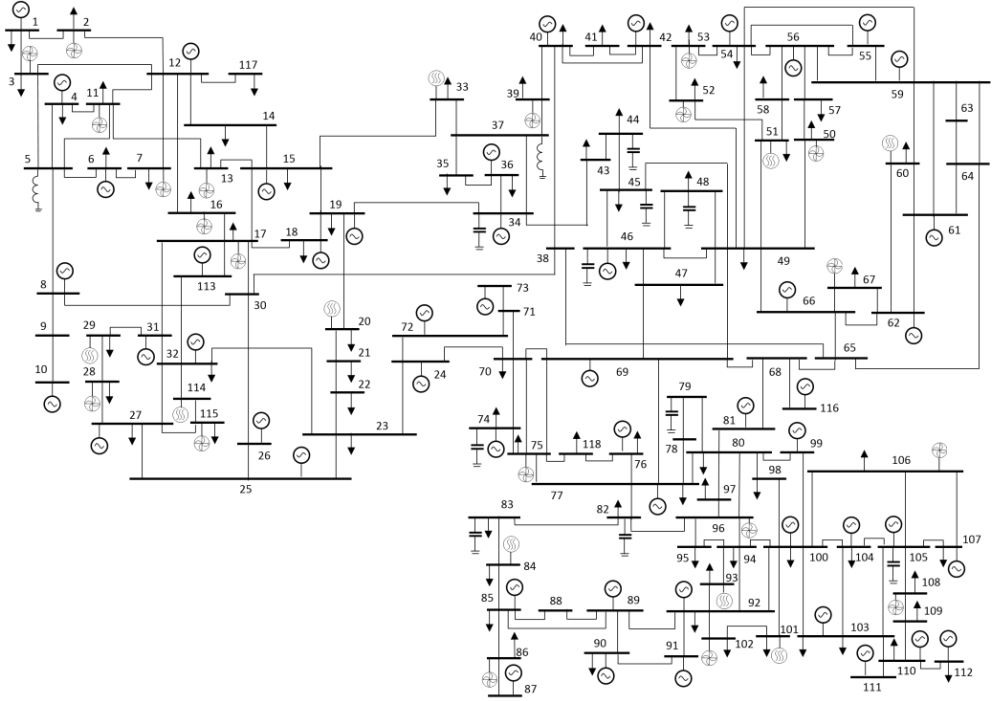

Figure 2: Modified IEEE 118-bus system diagram.

### 6.1.2 Linear Time Intervals for Short-term Scheduling

Short-term scheduling in this simulation is assumed to be the day-ahead scheduling. Typically, day-ahead generation scheduling is made per hour, therefore, the length of linear time intervals in our simulation is also chosen to be 1 hour within the entire 24-hour time period.

### 6.1.3 Linear Time Intervals for Real-time Dispatch

For the real-time dispatch, it is necessary to determine an appropriate timescale. The selection of the real-time dispatch timescale is closely related to factors such as the IPS power uncertainty, the IPS proportion, and the automatic generation control (AGC) power adjustment capability. Wang and Guo (2018) proposed an analytical formulation to determine the proper real-time dispatch critical timescale in power systems with different IPS proportions. For the power grid shown in Figure 2, according to the approach in (Wang and Guo 2018), the real-time dispatch timescale is set as 1.5 hours. Such a timescale is then divided into six linear time intervals of 15 minutes.

### 6.2  POPF Results for Short-term Scheduling

For the power grid in Figure 2, the TOPF and POPF models are applied into the short-term scheduling. We record the respective actual power generation costs as well as the numbers and magnitudes of power and voltage limits violations. For all the 24 time intervals, all operational constraints are satisfied in the POPF. The TOPF, on the contrary, produces notable violations of power/voltage magnitude over-limits, as depicted respectively in Figure 3(a) and Figure 3(b), wherein the "CVN" and "CVA" represent constraint violation numbers and amounts, respectively. It can be observed in Figure 3 that unlike the POPF, the TOPF results in power constraint violations in all the 24 time intervals and a number of nodal voltage constraint violations.

The actual generation costs of POPF and TOPF are compared in Figure 4. We can see that in general, POPF achieves lower generation cost than TOPF. In addition, there are over-limits in TOPF but POPF satisfies all constraints in this simulation.

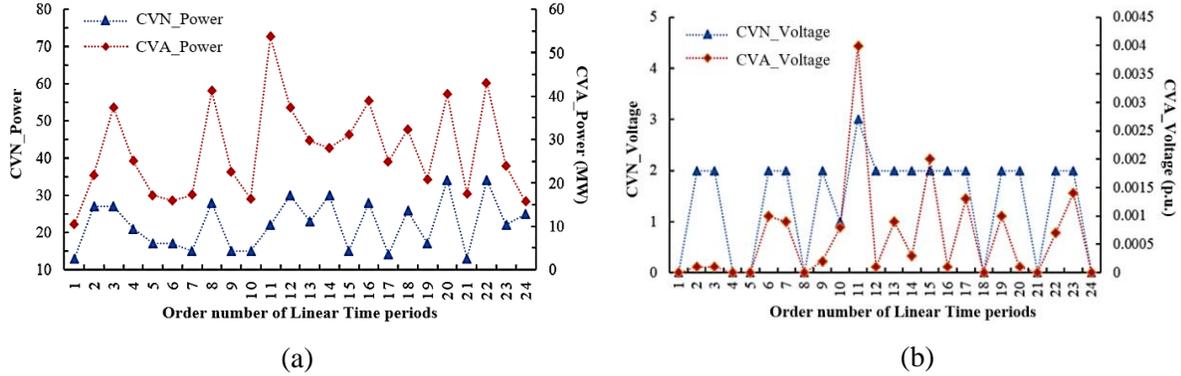

(a)          (b)

Figure 3: Over-limit number and over-limit amount of TOPF under day-ahead scheduling for (a) branch power, and (b) node voltages.

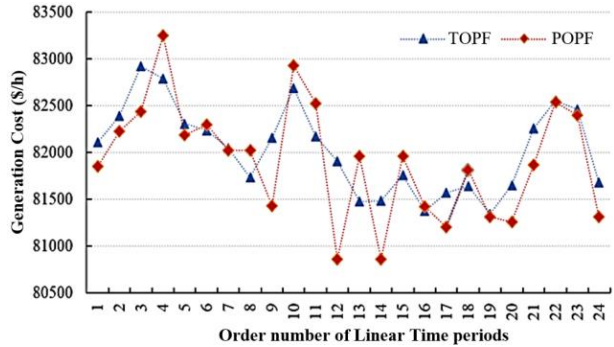

Figure 4: Generation Cost of TOPF and POPF under day-ahead scheduling.

### 6.3 POPF Results for Real-time Dispatch

Based on the connection between real-time dispatch and short-term scheduling, the TOPF and POPF models are also applied into real-time dispatch, the power generation cost and the constraint violations are thus obtained, respectively, with the results shown in Table 2. Note that since the real-time dispatch timescale is a rolling window, the initial time is selected to correspond to the beginning of the 10$^{th}$ time interval for short-term scheduling.

It can be observed from Table 2 that the POPF solution satisfies the constraints in all the six time intervals. The TOPF produces over-limits in both branch power flows and nodal voltage magnitudes. However, in comparison with the optimization results for short-term scheduling, the violation levels of TOPF are generally lower.

### 6.4 Discussion

The results above indicate that heuristically choosing discrete points may lead TOPF solution to higher generation cost and unpredictable over-limits in branch power flows and nodal voltage magnitudes. The proposed POPF, on the contrary, effectively alleviates these problems and produces more efficient generation schedules. This is because the POPF adaptively partitions the entire time period into several linear time intervals and takes median and two end snapshots of each linear time interval to ensure its performance over the entire time period.

In terms of the computational cost, POPF requires heavier computational cost than TOPF to investigate the optimality and feasibility based on the three representative snapshots. However, the raise of computational burden in the POPF is not significant compared to TOPF since only two additional snapshots are considered, whereas substantial gains in terms of optimality and security are finally achieved. In summary, the POPF method, which is convenient, effective, and efficient, is a more desirable method for power flow operations and planning over a time period.

Table 2: Simulation results of POPF and TOPF models in modified IEEE 118-bus system under real-time dispatch.

| Linear Time Intervals | Parameters | TOPF | POPF |
|---|---|---|---|
| 0h-0.25h | Power violated number(s)/Amount(MW) | 9/2.3479 | 0/0 |
| | Voltage violated number(s)/Amount(MW) | 2/0.0001 | 0/0 |
| | Generation cost($/h) | 80117 | 79953 |
| 0.25h-0.5h | Power violated number(s)/Amount(MW) | 14/4.0001 | 0/0 |
| | Voltage violated number(s)/Amount(MW) | 4/0.0004 | 0/0 |
| | Generation cost($/h) | 79852 | 79928 |
| 0.5h-0.75h | Power violated number(s)/Amount(MW) | 14/6.4342 | 0/0 |
| | Voltage violated number(s)/Amount(MW) | 2/0.0001 | 0/0 |
| | Generation cost($/h) | 79901 | 80049 |
| 0.75h-1h | Power violated number(s)/Amount(MW) | 15/4.3408 | 0/0 |
| | Voltage violated number(s)/Amount(MW) | 4/0.0001 | 0/0 |
| | Generation cost($/h) | 80096 | 80009 |
| 1h-1.25h | Power violated number(s)/Amount(MW) | 7/1.756 | 0/0 |
| | Voltage violated number(s)/Amount(MW) | 0/0 | 0/0 |
| | Generation cost($/h) | 79821 | 79923 |
| 1.25h-1.5h | Power violated number(s)/Amount(MW) | 15/2.8245 | 0/0 |
| | Voltage violated number(s)/Amount(MW) | 4/0.0001 | 0/0 |
| | Generation cost($/h) | 79934 | 79858 |

## 7 CONCLUSIONS

We propose a period optimal power flow (POPF) model to enhance the connection between the discrete snapshots employed in the optimization and the entire time period which the optimization results will be applied to. An important aspect of the linearization and the discretization in this approach is that the entire time period is partitioned into linear time intervals, which satisfies the following properties:

1) Power injections and nodal voltages within the time interval are linear functions of time;
2) The voltage magnitudes of the controllable power sources are constant in each linear time interval;
3) The interval construction technique ensures the Jacobian to be approximated as a constant matrix.

In all, in each linear time interval, we take its median snapshot to construct the objective function in the POPF and check operational limits for two end snapshots. Simulation case studies on a modified IEEE 118-bus system have demonstrated the simplicity and effectiveness of POPF. Compared to traditional optimal power flow, the proposed POPF achieves better performances and fewer unpredictable over-limits over a given time period.

## REFERENCES

Capitanescu, F. and Wehenkel, L. 2013. "Experiments with the interior-point method for solving large scale optimal power flow problems," Electric Power Systems Research, vol. 95, pp. 276-283.

Carpentier, J. 1962. "Contribution to the economic dispatch problem," Bull. Soc. France Elect., vol. 8, pp. 431–437.

Gill, S., Kockar, I. and Ault G. W. 2014. "Dynamic Optimal Power Flow for Active Distribution Networks," IEEE Transactions on Power Systems, vol. 29, no. 1, pp. 121-131.

Gerstner, P., Schick, M. and Heuveline, V. 2016. "A Domain Decomposition Approach for Solving Dynamic Optimal Power Flow Problems in Parallel with Application to the German Transmission Grid".

Kourounis, D., Fuchs, A., and Schenk, O. 2018. "Toward the Next Generation of Multiperiod Optimal Power Flow Solvers," IEEE Transactions on Power Systems, vol. 33, no. 4, pp. 4005-4014.

Liang, R., et al. 2011. "Optimal power flow by a fuzzy based hybrid particle swarm optimization approach," Electric Power Systems Research, vol. 81, no. 7, pp. 1466-1474.

Maria, G. A. and Findlay, J. A. 1987. "A Newton Optimal Power Flow Program for Ontario Hydro EMS," IEEE Transactions on Power Systems, vol. 2, no. 3, pp. 576-582.

Reddy, SS., Bijwe, PR. And Abhyankar, AR. 2015. "Real-Time Economic Dispatch Considering Renewable Power Generation Variability and Uncertainty Over Scheduling Period," IEEE Systems Journal, vol. 9, no. 4, pp. 1440-1451.

Ringlee, R. J. and Williams, D. D. 1962. "Economic System Operation Considering Valve Throttling Losses II-Distribution of System Loads by the Method of Dynamic Programming," Transactions of the American Institute of Electrical Engineers. Part III: Power Apparatus and Systems, vol. 81, no. 3, pp. 615-620.

Souza, RR., Balbo, AR., Nepomuceno, L., et al. 2017. "A primal-dual interior/exterior point method, with combined directions and quadratic test in reactive optimal power flow problems," IEEE Latin America Transactions, vol. 15, no. 8, pp. 1413-1421.

Sun, D. I., Ashley, B., Brewer, B., Hughes, A. and Tinney, W. F. 1984. "Optimal Power Flow by Newton Approach," IEEE Power Engineering Review, vol. PER-4, no. 10, pp. 39-39.

Trovato, V., Sanz, I. M., Chaudhuri, B. and Strbac, G. 2017. "Advanced Control of Thermostatic Loads for Rapid Frequency Response in Great Britain," IEEE Transactions on Power Systems, vol. 32, no. 3, pp. 2106-2117.

Wang, Z. and Guo, Z. 2017. "Toward a characteristic optimal power flow model for temporal constraints," 2017 IEEE Transportation Electrification Conference and Expo, Asia-Pacific (ITEC Asia-Pacific), Harbin, pp. 1-6.

Wang, Z. and Guo, Z. 2018. "On critical timescale of real-time power balancing in power systems with intermittent power sources," Electric Power Systems Research, vol. 155, pp. 246-253.

## AUTHOR BIOGRAPHIES

**Zongjie Wang** is a postdoctoral associate in department of Biological and Environmental Engineering at Cornell University. Her research interest include applications of simulation-optimizition techniques, renewable energy, power systems. Her email is zw337@cornell.edu.

**Ge Guo** is a postdoctoral associate in Systems Engineering at Cornell University. Her research interests include stochastic optimization with applications in sustainable power systems and manufacturing systems. Her email is gg442@cornell.edu.

**C. Lindsay Anderson** (M'01) is an Associate Professor, the Norman R. Scott Sesquicentennial Faculty Fellow, and the Kathy Dwyer Marble and Curt Marble Faculty Director for Energy with the Atkinson Center for a Sustainable Future, Cornell University, Ithaca, NY, USA. Her email is cla28@cornell.edu.